

\input amstex
\documentstyle{amsppt}

\magnification=1200
\pagewidth{5.4in}
\pageheight{7.5in}
\parskip 10pt

\expandafter\redefine\csname logo\string@\endcsname{}
\NoBlackBoxes                
\NoRunningHeads
\redefine\no{\noindent}

\define\C{\bold C}

\define\Z{\bold Z}

\redefine\P{\bold P}

\define\al{\alpha}
\define\be{\beta} 
\define\ga{\gamma}
\define\de{\delta}
\define\la{\lambda}
 
\define\Si{\Sigma} 
\define\La{\Lambda}

\define\Om{\Omega}

\define\th{\theta}
\define\om{\omega}

\define\na{\nabla}

\define\sub{\subseteq}

\define\st{\ \vert\ }   
\redefine\ll{\lq\lq}
\redefine\rr{\rq\rq\ }
\define\rrr{\rq\rq}

\redefine\deg{\operatorname {deg}}

\define\diag{\operatorname {diag}}

\redefine\dim{\operatorname {dim}}

\define\End{\operatorname {End}}

\define\Hom{\operatorname {Hom}}

\define\glnc{GL_n\C}
\define\Gl{GL}

\redefine\b{\partial}

\define\bbti{\frac{\b}{\b t_i}}
\define\bbqi{\frac{\b}{\b q_i}}

\topmatter

\title Quantum cohomology via D-modules
\endtitle

\author Martin A. Guest
\endauthor

\dedicatory
Dedicated to Graeme Segal on the occasion of his 60th birthday
\enddedicatory

\abstract
We propose a new point of view on quantum cohomology, motivated by 
the work of Givental and Dubrovin, but closer to differential geometry than the existing
approaches. The central object is a D-module which \ll quantizes\rr a commutative
algebra associated to the (uncompactified) space of rational curves. Under appropriate
conditions, we show that the associated flat connection may be gauged to the flat
connection underlying quantum cohomology. This method clarifies the role of the Birkhoff
factorization in the \ll mirror transformation\rrr, and it gives a new algorithm
(requiring construction of a Groebner basis and solution of a system of o.d.e.) for
computation of the quantum product. 
\endabstract

\endtopmatter

\document

Quantum cohomology
first arose in physics, and its (mathematically conjectural) properties
were supported by physical intuition.  A rigorous mathematical definition
came later, based on deep properties of certain moduli spaces.  
We shall propose another point of view on quantum cohomology, closer in spirit to
differential geometry.

The main ingredient in our approach
is a flat connection, considered as a holonomic D-module (or maximally overdetermined
system of p.d.e.). This object itself is not new:
Givental's \ll quantum cohomology D-module\rr is already well known (\cite{Gi1}),
and the associated flat connection appears in Dubrovin's 
theory of Frobenius manifolds (\cite{Du}).
But, in the existing literature, the D-module plays a subservient role, being a
consequence of the construction of the Gromov-Witten invariants
and the quantum cohomology algebra.  
For us, the D-module will be the main object of interest. 

We define a  quantization of a (commutative) algebra $\Cal A$
to be a (non-commutative) D-module $M^h$ which satisfies certain properties. 
The quantum cohomology D-module is
a particular kind of quantization, which arises in the
following way. For a K\"ahler manifold $M$, we start with an algebra $\Cal A$
which is associated to the \ll raw data\rr consisting of the set of all rational
curves in $M$.  Then we construct (or assume the existence of) a quantization $M^h$.
Next we transform $M^h$ into a new D-module $\hat M^h$ with certain properties.
Finally, de-quantization (\ll semi-classical limit\rrr)
produces a commutative algebra $\hat\Cal
A$, which (under appropriate conditions)
turns out to be the quantum cohomology $QH^\ast M$.

Our scope will be very modest in this article: we consider only the \ll small\rr 
quantum cohomology algebra
$QH^\ast M$ of a manifold $M$ whose ordinary cohomology algebra $H^\ast M$ is generated by
two-dimensional classes.  But this case is sufficiently nontrivial to demonstrate that
our method has something to offer, both conceptually and computationally.
The most obvious conceptual benefit is that 
the usual moduli space $\Cal M$  has been replaced by the D-module $M^h$. 
As a first application we give an algorithm for computing the structure
constants of the quantum cohomology algebra 
($3$-point genus zero Gromov-Witten invariants), in the case
of a Fano manifold. This involves a Gr\"obner basis calculation and
a finite number of \ll quadratures\rrr; it is quite different from
previously known methods.  A second application is a new interpretation of the
\ll mirror coordinate transformation\rrr. 
Impressively mysterious in its original context
(\cite{Gi3}--\cite{Gi4}, \cite{LLY1}--\cite{LLY3}), 
it arises here in a straightforward 
differential geometric fashion,   reminiscent of the well known transformation 
to local Euclidean coordinates for a flat Riemannian manifold.

Here is a more detailed description of the organization of this paper.  
In \S 1 we review some facts
concerning D-modules, mainly to establish notation.  In \S 2 we recall the
quantum cohomology algebra and the quantum product, 
again to set up notation.
\ll Quantum cohomology algebra\rr  refers to the isomorphism
type of the algebra, while \ll quantum product\rr means the product operation
on the vector space $H^\ast M$, i.e.\  a way of multiplying ordinary
cohomology classes.    

Our point of view is introduced in \S 3: we start with an algebra
$\Cal A$ and construct from it both a \ll quantum cohomology algebra\rr
and a \ll quantum product\rrr.
The method is conceptually straightforward. 
To a quantization $M^h$ of $\Cal A$ there corresponds
a flat connection $\nabla = d + \Om^h$, where $\Om^h$ has a simple pole
at $h=0$.  We may write $\Om^h = L^{-1}dL$ for some loop group-valued map $L$.
Replacing $L$ by $L_-$, where $L=L_-L_+$ is the Birkhoff factorization, we
obtain $\hat\Om^h = L_-^{-1}dL_- $, and the connection $d + \hat\Om^h$ is the
required connection.  The map $L$ is a (complicated) generating function for
certain Gromov-Witten invariants but we shall not need it.  Our main interest is the
gauge transformation $L_+=Q_0+O(h)$ which converts $\Om^h$ to $\hat\Om^h$. 
For the manifolds discussed here, $\Cal A$ and $M^h$ are known, and
$\Om^h$ can be computed.  If $L_+$ can be computed, then $\hat \Om^h$
(and the quantum cohomology algebra, together with its structure constants) can be
computed too.

In \S 4 we discuss the case of Fano manifolds.  Here it turns out that
$\Cal A = \hat\Cal A$, i.e.\  the \ll provisional\rr algebra is actually the
\ll correct answer\rrr.  The gauge transformation
$L_+$ has a special form but it is not trivial; 
indeed, its first term $Q_0$ tells us how to produce the
quantum product. Thus all quantum products can be {\it determined explicitly} 
by our method from the relations of the quantum cohomology algebra (more precisely,
from their quantizations). 
The following two families of manifolds are of special interest:

(1) Let $M=G/B$, the full flag manifold of a complex semisimple Lie group $G$.
The quantum cohomology algebra was found originally by Givental and Kim 
(\cite{Gi-Ki}, \cite{Ki}) and justified via the conventional moduli space theory. The first integrals
of the quantum Toda lattice provide a quantization $M^h$. It is known that the quantum
product can be described using quantum Schubert polynomials (see 
\cite{FGP}, \cite{Ki-Ma} for the case $G=\glnc$); therefore, the theory of such polynomials
is governed by our matrix $Q_0$. A more detailed treatment of flag manifolds from our point of
view can be found in \cite{Am-Gu}.

(2) Let $M$ be a Fano toric manifold.  In this
case a formula for the quantum cohomology was proposed by Batyrev (\cite{Ba1}), but
the subsequent proof of the correctness of the formula 
(see \cite{Co-Ka}, Chapter 11) depended on Givental's
mirror theorem from \cite{Gi4}.  The appropriate quantization is the
generalized hypergeometric D-module of \cite{GKZ} (whose relevance to
mirror symmetry was already known; cf.\  \cite{Ba2}, \cite{HLY}). 
Again, the matrix $Q_0$ produces the quantum product.

Beyond Fano manifolds there arises the interesting possibility that 
$\hat\Cal A$ may be different from $\Cal A$, and we discuss this in \S 5,
primarily with toric manifolds in mind.  Several authors
have pointed out that the quantum cohomology
algebra constructed by Batyrev in \cite{Ba1} is generally the
\ll wrong answer\rr  for a non-Fano toric manifold. Our point of
view resolves this apparent conflict,
at least in the case of semi-positive toric manifolds: Batyrev's algebra is $\Cal A$, the \ll usual\rr
quantum cohomology algebra is $\hat\Cal A$, and the two are
related via $L_+$. The gauge transformation $L_+$ contains more
information than in the Fano case, namely a coordinate transformation.
For toric manifolds this is Givental's mirror transformation.
It is a natural operation from the point of view of D-modules, 
but considerably less so from the point of view of the quantum cohomology algebra,
where it seems miraculous (\cite{Gi4}, \cite{Co-Ka}). 

The results of this paper can probably be generalized in
various directions.  For manifolds whose ordinary cohomology is not
generated by two-dimensional classes, one may work with the subalgebra
generated by such classes, as is standard in discussions of mirror
symmetry. For \ll big\rr quantum cohomology our methods may apply
to some extent.  Finally, there may well
be more general algebras $\Cal A$ to which our methods apply,
i.e.\  algebras without any obvious connection to quantum cohomology theory.

This project began with a conviction that
integrable systems methods could be used to
rehabilitate Batyrev's \ll incorrect\rr computations of quantum cohomology algebras
of toric varieties in \cite{Ba1}. 
It will be obvious to the experts that our framework owes much
to the ideas of Givental (\cite{Gi1}-\cite{Gi4}) and
Dubrovin (\cite{Du}), and we gladly acknowledge these as our
main sources of inspiration, though we would not have made much progress
without the excellent
treatments of quantum cohomology in \cite{Co-Ka} and hypergeometric
$D$-modules in \cite{SST}.

For background information on quantum cohomology we refer the reader to
the books \cite{Co-Ka},  \cite{Mn} and their references.  
In addition, survey articles related 
to the quantum differential equations include \cite{BCPP},  \cite{Pa},  \cite{Gu2}.
An introduction to loop group techniques in integrable systems can
be found in the book \cite{Gu1}.

The author is grateful to Josef Dorfmeister for suggesting
the use of the Lie algebra in the proof of Proposition 4.1, and to Hiroshi
Iritani for several very helpful comments.
He thanks Alexander Givental for explaining
that the present article has some overlap with the preprint of
T. Coates and A. Givental, \ll Quantum Riemann-Roch, Lefschetz and Serre\rrr,
math.AG/0110142.
The author was partially supported by a 
research grant from the JSPS.

\head
\S 1 D-modules and flat connections
\endhead

Let $K$ be an algebra of functions of the complex variables $q_1,\dots,q_r$.
(In practice we shall use the polynomial algebra $\C[q_1,\dots,q_r]$,
or the field of rational functions or germs of holomorphic functions.)
Depending on the
context, we regard $q_i$ either as a formal variable or as a function $t\mapsto q_i=e^{t_i}$
where $t=(t_1,\dots,t_r)\in\C^r$.  We introduce the notation $\b_i = \bbti= q_i\bbqi$, and
define $D$ to be the algebra of differential operators generated by 
$\frac{\b}{\b q_1},\dots,\frac{\b}{\b q_r}$ with
coefficients in $K$.
Let $M=D/(D_1,\dots,D_u)$ be a cyclic D-module (a left module over $D$, generated
by the constant differential operator $1$), where $(D_1,\dots,D_u)$ means the
left ideal generated by differential operators $D_1,\dots,D_u$.  In this section
we shall assume that
$M$ is free over $K$ of rank $s+1$.    For basic facts on
D-modules we refer to \cite{SST},  \cite{Ph},  \cite{Co}.

The D-module $M$ is an algebraic version of the system of partial differential
equations $D_1f=\dots=D_uf=0$.  Here, $f$ belongs to a given function space $F$,
but $M$ is of course independent of $F$ (and this is its advantage).  To say
that $M$ has finite rank over $K$ is to say, roughly speaking, 
that the system is \ll maximally overdetermined\rrr;
in particular its solution space is finite dimensional.  More precisely, the vector
space $\Hom_D(M,F)$ is called the {\it solution space of} $M$ with respect to
the function space $F$, and this is isomorphic
to the usual solution space $\{f\in F\st D_1f=\dots=D_uf=0\}$ of the system:  to a solution
$f$ there corresponds the D-module homomorphism $M\to F$ given by $P\mapsto Pf$
(for any $P\in D$).  The solution space (in either sense) is a complex vector space of
dimension $s+1$.

We shall review briefly the relation between D-modules and flat
connections.  
Let us choose differential operators $P_0,\dots,P_s$ such that the equivalence
classes $[P_0],\dots,[P_s]$ form a $K$-module basis of $M$. 
(There is a standard way of doing this, by constructing first a
Gr\"obner basis of the ideal $(D_1,\dots,D_u)$, as explained in Section 1.4 of \cite{SST}.)
Without loss of generality
we may assume $P_0=1$.  With respect to this basis we define matrices 
$\Om_i=(\Om^{i}_{kj})_{0\le k,j\le s}$ by
$$
[\b_i P_j] = \sum_{k=0}^s \Om^{i}_{kj} [P_k],
$$
and we put $\Om=\sum_{i=1}^r \Om_i dt_i$, a $1$-form with values in the
space $\End(\C^{s+1})$ of complex $(s+1)\times(s+1)$ matrices.  
The formula $\na=d+\Om$ defines a connection in the trivial vector bundle
$\C^r\times \C^{s+1} \to\C^r$, where $\C^{s+1}$ is identified with the vector space spanned by
$[P_0],\dots,[P_s]$.  Namely, $\na_{\b_i}[P_j] =
\sum_{k=0}^s \Om^{i}_{kj} [P_k]$,  
and more generally for any section $\sum_{j=0}^s y_j [P_j]$ of this bundle,
$
\na_{\b_i}(\sum_{j=0}^s y_j [P_j]) =
\sum_{j=0}^s \b_i y_j [P_j] + \sum_{j=0}^s y_j \na_{\b_i}[P_j].
$

\proclaim{Proposition 1.1} The connection $\na$ is flat.
\endproclaim

\demo{Proof} By definition we have $\na_{\b_i}\na_{\b_j}=\na_{\b_j}\na_{\b_i}$
(since $\b_i \b_j = \b_j \b_i$), so the curvature tensor of
$\na$ is zero.
Alternatively, the zero curvature condition
$d\Om + \Om\wedge\Om = 0$ follows directly from computing
both sides of the equation $\b_i \b_j [P_k] = \b_j \b_i [P_k]$.
\qed\enddemo

\proclaim{Proposition 1.2} We have an isomorphism of vector spaces
$$
\Hom_D(M,F)\longrightarrow \{ \text{covariant constant sections of $\na^\ast$} \},
\quad
f\longmapsto
\pmatrix
P_0 f\\
\vdots\\
P_s f
\endpmatrix
$$
where $\na^\ast$ is the dual connection to $\na$.
\endproclaim

\demo{Proof} On the left hand side, $f$ is regarded as the D-module
homomorphism $P\mapsto Pf$, whereas on the right hand side $f$ is
a solution of the system $D_1 f = \dots = D_u f = 0$. The dual connection
is defined by $(\na_{\b_i}^\ast [P_j]^\ast)[P_{k}] = - [P_j]^\ast (\na_{\b_i} [P_{k}])$
where $[P_0]^\ast,\dots,[P_s]^\ast$ is the dual basis to $[P_0],\dots,[P_s]$.
The column vector in the statement of the proposition refers to the
section $\sum_{j=0}^s (P_jf) [P_j]^\ast$.
A section $\sum_{j=0}^s y_j [P_j]^\ast$ is covariant constant if
the following expression is zero for all $k$:
$$
\align
(\na_{\b_i}^\ast \sum_{j=0}^s  y_j [P_j]^\ast)[P_{k}]
&= ( \sum_{j=0}^s {\b_i} y_j [P_j]^\ast
+  \sum_{j=0}^s  y_j \na_{\b_i}^\ast [P_j]^\ast) [P_{k}]\\
&= \b_i y_{k}  -   
\sum_{j=0}^s  y_j ( [P_j]^\ast \sum_{l=0}^s \Om^i_{lk}[P_{l}])\\
&= \b_i y_{k}  - \sum_{j=0}^s  y_j \Om^i_{jk}.
\endalign
$$
For any $f\in \Hom_D(M,F)$, we have to verify that
$y_{k}= P_{k} f$ defines a covariant constant section.  But this follows
immediately from the formula $[\b_i P_{k}] = \sum_{k=0}^s \Om^{i}_{jk} [P_j]$
defining $\Om$.  The map in question is therefore a well defined, linear, map.
To prove that it is an isomorphism, we observe that the kernel is zero
(because $P_0 f=f$), and that
$\dim \Hom_D(M,F) = s+1$ by assumption.
\qed \enddemo

This generalizes the well known elementary construction of a system of
first order o.d.e.\  equivalent to a higher order o.d.e.  Here we construct the system
$\b_i y_{k}  = \sum_{j=0}^s  y_j \Om^i_{jk}$ of first order p.d.e.\  equivalent to
the higher order system $D_1f=\dots=D_uf=0$. 
Conversely, given a flat connection (hence a system of
first order p.d.e.), it is possible to construct a cyclic 
D-module of finite rank over an appropriate algebra $K$
(hence a system of higher order p.d.e.).

Since the dual connection $\nabla^\ast = d-\Om^t$ is flat, there exist covariant constant
sections $H_0,\dots,H_s$ which are linearly independent at each point of $\C^r$.
Representing these sections by column vectors, as above, let us introduce
$$
H=
\pmatrix
\vert &  & \vert \\
H_0 & \cdots & H_s \\
\vert &  & \vert
\endpmatrix
$$
i.e.\  the \ll fundamental solution matrix\rr of the first order system.  By definition
we have $\Om^t = dH H^{-1}$.  Up to multiplication on the right by a constant
invertible matrix, this equation determines $H$ uniquely.
Equivalently, if $f_0,\dots,f_s$ are a basis of solutions of the higher order system
$D_1f=\dots=D_uf=0$, and if $J=(f_0,\dots,f_s)$ is regarded as
a row vector, then
$$
H=
\pmatrix
- & P_0 J & - \\
 & \vdots &  \\
- & P_s J & -
\endpmatrix
$$
satisfies  $\Om^t = dH H^{-1}$.

A standard technique is to study the transformation (symbol map)
$\b_i\mapsto b_i$
from the
non-commutative algebra $D$ to the commutative algebra $K$.
A differential operator $P$ maps to a polynomial $\tilde P$.  The D-module
$M=D/(D_1,\dots,D_u)$ is transformed to a $K$-module
$\tilde M=K[b_1,\dots,b_r]/(\tilde D_1,\dots,\tilde D_u)$, and 
the associated flat connection $\na$ is transformed to a
connection $\tilde \na$, {\it but the connection}
 $\tilde \na$ {\it  is not in general flat.}
In more detail, we have $\na=d+\sum_{i=1}^r\Om_i dt_i$ where $\Om_i$
is the matrix representing the action of the differential operator $\b_i$, and $\tilde\na = d
+ \sum_{i=0}^r \tilde\Om_i dt_i$ where $\tilde\Om_i$ is the matrix representing 
the action of the operator $b_i$.   
As explained earlier, the fact that $\b_i \b_j = \b_j \b_i$ 
leads to the flatness condition $d\Om + \Om\wedge\Om = 0$.  However,
the condition $b_i b_j = b_j b_i$ says only that $\tilde\Om\wedge \tilde\Om = 0$.
The exterior derivative $d\tilde\Om$ is not in general zero.  This phenomenon
is the key to our construction of quantum cohomology in \S 3.

\head
\S 2 The quantum cohomology D-module
\endhead

In this section we shall review briefly the Dubrovin connection (or D-module)
which arises in the standard construction of quantum cohomology theory.
We begin with a compact K\"ahler manifold $M$ of (complex) dimension
$n$, whose ordinary cohomology algebra --- with complex coefficients ---
is of the form
$$
H^\ast M = \C[b_1,\dots,b_r]/(R_1,\dots,R_u)
$$
where $b_1,\dots,b_r$ are additive generators of $H^2 M$ and
$R_1,\dots,R_u$ are certain relations (polynomials in $b_1,\dots,b_r$).
(As mentioned in the introduction, this assumption can be removed
by studying the subalgebra generated by two-dimensional cohomology
classes.)
By general principles it follows that the (small) quantum cohomology
algebra is of the form
$$
QH^\ast M = K[b_1,\dots,b_r]/(\Cal R_1,\dots,\Cal R_u)
$$
where $K=\C[q_1,\dots,q_r]$ and each $\Cal R_i$ is a \ll $q$-deformation\rr
of $R_i$.  (For certain $M$, an extension or completion
of $K$ may be necessary here, but we shall assume in this
section that $M$ is not of this type.)  
As in \S 1, the variables $q_1,\dots,q_r$ here may be
considered either as formal variables or as functions
$
q_i:t= \sum_{j=1}^r t_j b_j \mapsto  e^{t_i}
$
on $H^2M$.
With the latter convention,
$H^\ast M$ and $QH^\ast M$ are isomorphic
as vector spaces (but not, in general, as algebras), for each value of $t$.

Quantum cohomology theory gives, in addition to $QH^\ast M$,
a quantum product operation on $H^\ast M$. That is, for any $x,y\in H^\ast M$, there is 
an element 
$x\circ_t y\in H^\ast M$, which has the property
$x\circ_t y=x\cdot y + \text{terms involving $q_i=e^{t_i}$, $1\le i\le r$}$,
where $x\cdot y$ denotes the cup product.  The relations $R_1,\dots,R_u$
are those of the algebra $(H^\ast M, \ \cdot\ )$, while the 
relations $\Cal R_1,\dots,\Cal R_u$
are those of the algebra $(H^\ast M, \ \circ_t\ )$.
In particular this gives rise to an isomorphism of vector spaces
$\de:QH^\ast M\to H^\ast M$ which \ll evaluates\rr a polynomial using
the quantum product.

The Dubrovin connection is the (complex) connection $\na=d +\frac 1h \om$
on the trivial bundle  $\C^r\times \C^{s+1} \to\C^r$ where $\om$ is
the complex $\End\,\C^{s+1}$-valued $1$-form on $\C^r$ defined by
$\om_t(x)(y) = x\circ_t y$.  Here $h$ is a nonzero complex parameter, so
in fact we have a family of connections.

\proclaim{Theorem 2.1} For any $h$ the connection $\na=d +\frac 1h \om$ is flat,
i.e.\  $d\om = \om\wedge\om = 0$.
\qed
\endproclaim

\no A proof of this well known theorem and further explanation can be found in \cite{Co-Ka} and
the other references on quantum cohomology at the end of this paper. 

\head
\S 3 Reconstructing quantum cohomology
\endhead

We begin with an abstract algebra of the form
$$
\Cal A = K[b_1,\dots,b_r]/(\Cal R_1,\dots,\Cal R_u),
$$
where the relations $\Cal R_1,\dots,\Cal R_u$ are homogeneous
with respect to a fixed assignment of degrees $\vert b_i\vert$,
$\vert q_j\vert$.  We shall always choose 
$\vert b_1\vert=\dots=\vert b_r\vert=2$, but 
$\vert q_1\vert,\dots,\vert q_r\vert$ (not necessarily non-negative)
will be specified later.  In
addition we assume that $\Cal A$  is a free $K$-module of rank $s+1$.
Finally, we assume that $\Cal A$ is a deformation of an algebra
$$
\Cal A_0  = \C[b_1,\dots,b_r]/(R_1,\dots,R_u)
$$
in the sense that $\Cal R_i\vert_{q=0} = R_i$
for $i=1,\dots,u$ and $\dim_{\C} \Cal A_0 = s+1$. 

Although it will play no role in this section, we should mention
that the situation we have
in mind is where $\Cal A_0=H^\ast M$ for
a compact connected K\"ahler manifold $M$, and where 
$\Cal A$ is obtained by using the (uncompactified) space of
rational curves in $M$ to define 
structure constants in the \ll naive\rr way as in early papers
in the physics literature.
In our examples $M$ will be a flag manifold $G/B$ or a toric manifold,
and we shall specify $\Cal A$ precisely when we discuss those cases.

Our main objective in this section will be to construct connections
satisfying the property of Theorem 2.1.  For this purpose,
we introduce the ring $D^h$ of differential operators 
generated by $h\b_1,\dots,h\b_r$ with coefficients in $K[h]$, and we
make the following fundamental definition:

\proclaim{Definition 3.1} A {\it quantization} of $\Cal A$ is a
D-module $M^h= D^h/(D_1^h,\dots,D_u^h)$ such that

\no (1) $M^h$ is free over $K[h]$ of rank $s+1$, 

\no (2) $\lim_{h\to 0} S(D^h_i) = \Cal R_i$,
where $S(D^h_i)$ is the result of replacing $h\b_1,\dots,h\b_r$
by $b_1,\dots,b_r$ in $D^h_i$ (for $i=1,\dots,u$).
\endproclaim

\no This notion depends on the specified generators and relations of $\Cal A$,
of course.
There is no guarantee
that such a quantization exists, but
it is sometimes possible to produce a quantization simply by replacing 
$b_1,\dots,b_r$ by $h\b_1,\dots,h\b_r$ in each $\Cal R_i$. When this works,
i.e.\  when the resulting
D-module is free of rank $s+1$, we refer to it as 
the {\it naive quantization}.

Assume now that $M^h$ is a quantization of $\Cal A$.
Then we may choose a $K[h]$-module basis
$[P_0],\dots,[P_s]$ of
$M^h$ such that $[c_0=\lim_{h\to 0} S(P_0)], \dots, [c_s=\lim_{h\to 0} S(P_s)]$
is a $K$-module basis of $\Cal A$.  We shall always do this by taking
$P_0,\dots,P_s$ to be the \ll standard monomials\rr 
in $h\b_1,\dots,h\b_r$ with respect to a choice
of Gr\"obner basis for the ideal $(D_1^h,\dots,D_u^h)$. 
For definiteness we
use the graded reverse lexicographic monomial order in which
$\b_1,\dots,\b_r$ are assigned weight one with
$\b_1>\dots>\b_r$.
(Gr\"obner basis theory for this situation is
explained in \cite{SST}.  Explicit computations may be carried
out using the Ore algebra package of the software Maple, \cite{Ma}.)
We define a connection form $\Om^h=\sum_{i=1}^r \Om^h_i dt_i$
as follows:

\proclaim{Notation} For $i=1,\dots,r$:

\no(1) let $\Om^h_i$
denote \ll the matrix of the action of $\b_i$\rr on the $K[h]$-module $M^h$,
i.e.\  $[\b_i P_j]=\sum_{k=0}^s (\Om^h_i)_{kj}[P_k]$;

\no(2) let $\om_i$
denote the matrix of multiplication by $b_i$ on the $K$-module
$\Cal A$, i.e.\  $[b_ic_j]=\sum_{k=0}^s (\om_i)_{kj}[c_k]$.
\endproclaim

\no It follows that 
$h\Om^h$ is polynomial in $h$, so $\Om^h$ is of the form
$$
\Om^h = \frac 1h \om + \th^{(0)} + h\th^{(1)} + \dots + h^p \th^{(p)},
$$
where $\om=\sum_{i=1}^r \om_i dt_i$,
$\th^{(0)},\dots,\th^{(p)}$ are matrix-valued $1$-forms,   and
$p$ is a non-negative integer which depends on the relations
$\Cal R_1,\dots, \Cal R_u$. 

If $\th^{(0)},\dots,\th^{(p)}$ were all zero, then 
the connection $\na=d+\Om^h$ (which is flat, by \S 1) would satisfy the condition of
Theorem 2.1, and hence would be a candidate for the Dubrovin connection.
It turns out that this situation can be achieved by making a suitable modification:

\proclaim{Proposition 3.2} Assume that $\Om^h$ depends holomorphically
on $q=(q_1,\dots,q_r)$, for $q$ in some open subset $V$. 
Then, for any point $q_0$ in $V$, there is a neighbourhood $U_0$ of $q_0$
on which the connection $\na = d + \Om^h$ is 
gauge equivalent to a connection $\hat\nabla =d + \hat\Om^h$
with $\hat\Om^h = \frac 1h \hat\om$, $\hat\om= Q_0 \om Q_0^{-1}$,
for some holomorphic map $Q_0:U_0\to \Gl(\C^{s+1})$.
\endproclaim

\demo{Proof} Since $d + \Om^h$ is flat, we have $\Om^h = L^{-1}dL$
for some $L:V\to \La \Gl(\C^{s+1})$.  
(In the notation of \S 1, $L=H^t$.)
Here, $\La \Gl(\C^{s+1})$ is the
(smooth) loop group of $\Gl(\C^{s+1})$, i.e.\  the space of all (smooth)
maps $S^1\to \Gl(\C^{s+1})$, where $S^1 = \{ h\in\C \st \vert h\vert = 1 \}$.
Let $L=L_-L_+$ be the Birkhoff factorization of $L$, where $L_+$ extends 
holomorphically to the disc $0\le \vert h\vert < 1$ and $L_-$ to
the disc $1<\vert h\vert \le\infty$, and where $L_-\vert_{h=\infty}=I$.
This factorization exists if and only if $L$ takes
values in the \ll big cell\rr of the loop group. For any given point 
$q_0$ of $V$, we may choose $\ga\in \La \Gl(\C^{s+1})$ so
that $\ga L(q_0)$ belongs to this big cell.  Replacing $L$
by $\ga L$, we obtain a factorization at $q_0$, and hence
on a neighbourhood $U_0$ of this point.
We may write
$$
\align
L_-(q,h)&=I+h^{-1}A_1(q)+h^{-2}A_2(q)+\dots\\
L_+(q,h)&= Q_0(q)(I + h Q_1(q)  + h^2 Q_2(q)  + \dots)
\endalign
$$
for some 
$A_i,Q_j:U_0\to \Gl(\C^{s+1})$.

Now we employ a well known argument from the theory of integrable
systems.   The gauge transformation
$L\mapsto \hat L = L(L_+)^{-1}=L_-$ transforms $\Om^h = L^{-1}dL$
into $\hat\Om^h = \hat L^{-1}d\hat L = L_-^{-1} dL_{-}$, and
the Laurent expansion of the latter manifestly contains only
negative powers of $h$.  But we have the alternative expression
$$
\align
L_-^{-1} dL_{-} &= (LL_+^{-1})^{-1} d(L L_+^{-1})
= L_+ L^{-1} dL L_+^{-1}  +  L_+ d(L_+^{-1}) \\
&= L_+ (\frac 1h \om + \th^{(0)} + h\th^{(1)} + \dots + h^p \th^{(p)})
L_+^{-1}  +  L_+ (dL_+^{-1}),
\endalign
$$
whose only negative power of $h$ occurs in the term $\frac 1h Q_0 \om Q_0^{-1}$.  
It follows that $\hat \Om^h = \frac 1h Q_0 \om Q_0^{-1}$, as required.
\qed\enddemo

Another way to express this modification is to say that we replace the original
basis $[P_0],\dots,[P_s]$ of $M^h$ by a new basis
$[\hat P_0],\dots,[\hat P_s]$, where
$\hat P_i = \sum_{j=0}^s (L_+)^{-1}_{ji} \hat P_j$.  Then $\hat \Om^h_i$ is
the matrix of the action of $\b_i$ with respect to the basis
$[\hat P_0],\dots,[\hat P_s]$.
At the same time, we replace the original basis $[c_0],\dots,[c_s]$ of
$\Cal A$ by the new basis $[\hat c_0],\dots,[\hat c_s]$, where
$\hat c_i = \sum_{j=0}^s (Q_0^{-1})_{ji} c_j$;  $\hat\om_i$ is
the matrix of multiplication by $[b_i]$  with respect to this new basis.
In this description, the entries of $(L_+)^{-1}$ are assumed to lie in $K[h]$.

The modified connection $\hat\na = d + \hat \Om^h$ will be the basic ingredient in
our construction of a \ll quantum cohomology algebra\rr $\hat\Cal A$
and a \ll quantum product operation\rrr.  The construction will be given here
in a special case, the general case being postponed to \S 5.  Namely,
we assume that 
$$
\text{$\hat c_0 = c_0 = 1$ and
$\hat c_i = c_i = b_i$ for $1\le i\le r$,}
$$
and that $L_+\vert_{q=0}=I$ ($L_+$ is then determined uniquely).
In this situation we simply define $\hat\Cal A=\Cal A$.  The 
\ll quantum product operation\rr will be defined on $\Cal A_0$, and
for this it is convenient to introduce the following terminology.

\proclaim{Notation}  For a polynomial $c$ in $b_1,\dots,b_r,q_1,\dots,q_r$
we denote the corresponding element of $\Cal A$ ---
the equivalence class of $c$ mod $\Cal R_1,\dots,\Cal R_u$ --- by $[c]$.
If $c$ is a polynomial in $b_1,\dots,b_r$
we denote the corresponding element of $\Cal A_0$ by $[[c]]$.
\endproclaim

\no We define
$$
\de:\Cal A\to \Cal A_0,\quad
[\hat c_i] \mapsto [[\hat c_i\vert_{q=0}]]\quad
(0\le i\le r).
$$
This is obviously an isomorphism of vector spaces if $q_1,\dots,q_r$
are considered as functions (and if $q_1,\dots,q_r$ are considered as formal
variables, $\de$ defines an isomorphism of $K$-modules 
$\Cal A\to \Cal A_0\otimes K$).  We introduce a \ll quantum product 
operation\rr $\circ_t$ on $\Cal A_0$ as follows:
$$
x\circ_t y= \de(\de^{-1}(x) \de^{-1}(y) ).
$$
(For a discussion of the relation between $\de$ and $\circ_t$,  see \S1 of \cite{Am-Gu}.)    
It follows that  the matrix of the operator $b_i\circ_t$ on $\Cal A$,
with respect to the basis $[\hat c_0],\dots,[\hat c_s]$,
is $\hat\om_i$, and hence that the \ll Dubrovin connection\rr associated
to $\circ_t$ is $d+\frac 1h \hat\om$. This is flat (since the gauge
equivalent connection $d+\Om^h$ is
flat, by \S 1), and so it satisfies $d\hat\om = \hat\om\wedge\hat\om=0$.

We postpone to later sections a discussion of when our abstract quantum
product coincides with the usual quantum product.  For the moment we wish to
emphasize that we have constructed a product with the expected properties,
and that our construction involves {\it a priori} the following steps: (1) an algebraic
(Gr\"obner basis) calculation to find $\Om^h$; (2) solution of a system of
ordinary differential equations to find $L$; (3) the factorization
$L=L_-L_+$.   Although steps (2) and (3) seem formidable in general, we shall see that they	
can sometimes  be reduced to a straightforward algorithm.

We conclude this section by giving some general properties of $\Om^h$.
Let $M^h_i$ be the subspace of $M^h$ which is spanned (over $K[h]$)
by the 
basis vectors $P_j$ of degree $i$ in $h\b_1,\dots,h\b_r$. Then
we have a decomposition
$
M^h  = M^h_0 \oplus M^h_1  \oplus \dots \oplus M^h_v,
$
with respect to which the $(\al,\be)$-th block of the matrix $\Om^h_i$
will be denoted $(\Om^h_i)_{\al,\be}$.  We
shall generally  use Greek indices, separated by commas, in reference
to block matrices.

\proclaim{Proposition 3.3} (1) For $\al\ge\be+2$ we have $(\Om^h)_{\al,\be}=0$.

\no Assume that the generators $D^h_i$ are homogeneous
in $h$, $q_1,\dots,q_r$, $\b_1,\dots,\b_r$, where: 
$h$ is assigned degree $2$, 
$q_1,\dots,q_r$ have their usual degrees, and $\b_1,\dots,\b_r$
are assigned degree $0$.  Then:

\no(2a) Each nonzero entry of the block $(\Om^h_i)_{\al,\be}$ has degree $2(\be-\al)$.

\no Assume further that $L_+\vert_{q=0}=I$.  Then:

\no(2b) Each nonzero entry of the block 
$(L_+)_{\al,\be}$ has degree $2(\be-\al)$.  In particular
each nonzero entry of
$(Q_i)_{\al,\be}$ has degree $2(\be-\al-i)$.
\endproclaim

\demo{Proof} (1) It follows from the division algorithm that
the filtration of $M^h$ defined by
$M^h_{(j)} = \oplus_{k=0}^j M^h_k$  satisfies
$h\b_i M^h_{(j)} \sub M^h_{(j+1)}$.
(2a) This is immediate from the definition of $\Om^h$
and the homogeneity of the $D^h_i$.
(2b) The homogeneity property of $\Om^h$ can be expressed as
$$
\Om^h(q_1,\dots,q_r)=
\diag(\la^{2v},\la^{2v-2},\dots,1)
^{-1}
\Om^{\la^2h}(\la^{\vert q_1\vert}q_1,\dots,\la^{\vert q_r\vert}q_r)
\diag(\la^{2v},\la^{2v-2},\dots,1)
$$
where $\diag(\la^{2v},\la^{2v-2},\dots,1)$ denotes a matrix
in block diagonal form.
We must show that the function $L_+$ satisfies the same condition.  
By the proof of Proposition 3.2, $L_+$ is determined uniquely by the differential
equation
$$
\frac1h Q_0\om Q_0^{-1} L_+ \ =\ L_+ \Om^h \ -\ dL_+
$$
and the condition $L_+\vert_{q=0}=I$.  Therefore, it suffices to observe that 
$$
\diag(\la^{2v},\la^{2v-2},\dots,1)
^{-1}
L_+(\la^{\vert q_1\vert}q_1,\dots,\la^{\vert q_r\vert}q_r,\la^2 h)
\diag(\la^{2v},\la^{2v-2},\dots,1)
$$
satisfies the same conditions.  
\qed\enddemo

\head
\S 4 Fano manifolds
\endhead

It is well known that
a Fano manifold, by which we mean a K\"ahler manifold $M$ whose K\"ahler
$2$-form represents the first Chern class $c_1 M$ of the manifold, has
particularly well behaved quantum cohomology.  
It is natural to begin by applying the theory of \S 3 in this case.

We start with a deformation
$\Cal A = K[b_1,\dots,b_r]/(\Cal R_1,\dots,\Cal R_u)$ of the
cohomology algebra $\Cal A_0 = H^\ast M = \C[b_1,\dots,b_r]/(R_1,\dots,R_u)$.
(A priori, $\Cal A$ may or may not be isomorphic to
the quantum cohomology algebra.)  For $G/B$ and toric manifolds,
suitable algebras $\Cal A$, and, most importantly, their quantizations
$M^h$, are already available \ll off the shelf\rrr.
Before looking at these in more detail, we shall point out some
further properties of the connection form $\Om^h$ in the Fano case.  
A basic ingredient is the fact that, from
the naive construction of $\Cal A$ using rational curves, 
the degree of $q_i$ satisfies
$
\vert q_i\vert \ge 2.
$
In the case of flag manifolds and Fano toric manifolds, 
this property leads to operators $D^h_i$ of the form
$h^{\vert I\vert}\b_{I} +$ lower order terms, where $\vert I\vert \ge 2$
and the lower order terms
have coefficients in the polynomial algebra $K[h]=\C[q_1,\dots,q_r,h]$;
we shall say that such $D^h_i$ are \ll regular\rrr.
It follows from this and the homogeneity property that the elements of
the Gr\"obner basis are also regular, and hence that $M^h$ is free over $K[h]$.
The matrices $h\Om^h_i$ will then have entries in $K[h]$.

\proclaim{Proposition 4.1} 
Assume that  $\Om^h_1,\dots,\Om^h_r$ are polynomial in $q_1,\dots,q_r$
with  $\vert q_1\vert,\dots,\vert q_r\vert \ge 4$. Then 
$L_+=Q_0(I + h Q_1  + h^2 Q_2  + \dots)$ satisfies:

\no(1) $Q_0=\exp X$ where $X_{\al,\be}=0$ for $\al\ge\be-1$,

\no(2) for $i\ge 1$, $(Q_i)_{\al,\be}=0$ for $\al\ge\be-i-1$.

\no In particular,  $Q_i=0$ for $i$ sufficiently large, i.e.\ 
$L_+$ must be a polynomial in $h$.
\endproclaim

\demo{Proof}   Since 
$h\Om^h =  \om + h\th^{(0)} + h^2\th^{(1)} + \dots + h^{p+1}\th^{(p)}$,
it follows from the homogeneity and polynomiality properties that $\theta_i^{(j)}$ satisfies 
$(\theta_i^{(j)})_{\al,\be}=0$
for $\al\ge \be-j-1$. 
Hence $\Om^h$ takes values in the Lie algebra
consisting of loops of the form 
$\sum_{i\in\Z}h^i A_i$ whose coefficients satisfy the following conditions: 
$(A_i)_{\al,\be}=0$ for 
$\al\ge \be-i-1$ when $i\ge 0$, and
$(A_i)_{\al,\be}=0$ for $\al\ge \be-i+1$ when $i< 0$.
Hence $L$ and $L_-,L_+$ take values
in the corresponding loop group.  In particular
$(L_+)^{-1}dL_+=\sum_{i\ge 0}h^i A_i$ where
$(A_i)_{\al,\be}=0$ for 
$\al\ge \be-i-1$,
from which the stated properties of $L_+$ follow.  
\qed\enddemo

\proclaim{Corollary 4.2} With the assumptions of Proposition 4.1,
we may assume that 
$\hat c_0 = c_0 = 1$ and $\hat c_i = c_i=b_i$ for $i=1,\dots,r$.
\endproclaim

\demo{Proof}  
We can assume
that $P_0=1$ and $P_i=h\b_i$ for $1\le i\le r$ (as a nontrivial relation
between  $h\b_1,\dots,h\b_r$ would lead to a nontrivial relation
between $b_1,\dots,b_r$). 
Hence we may take
$c_0=1$ and $c_i=b_i$ for $1\le i\le r$.  Next, by (1) of Proposition 4.1,
we have
$$
Q_0=
\pmatrix
1 & 0 & [\ast] \\
0 & I & [\ast] \\
[0] & [0] & [\ast]
\endpmatrix
$$
where  $[\ast]$ denotes a submatrix and
$[0]$ denotes a zero submatrix (where a submatrix may
consist of several blocks). Thus, $\hat c_i = c_i$ for $i=0,\dots,r$.
\qed\enddemo

This means that we are in the situation of \S 3:  we can define $\hat\Cal A=\Cal A$ and
we obtain a \ll quantum product operation\rr on $H^\ast M$ from $L_+$.
If $L_+\vert_{q=0}=I$, $L_+$ is homogeneous, by (2b) of Proposition 3.3.
As each $\vert q_i\vert$ is positive, homogeneity implies that $L_+$ is polynomial
in each $q_i$, so the procedure of \S 3 gives a change of basis of $\Cal A$.
Moreover, Proposition 4.1 leads to an explicit algorithm
for $L_+$.    The essential point is
that $L_+$ is characterized by the system of equations
$$
\frac1h Q_0\om Q_0^{-1} L_+ \ =\ L_+ \Om^h \ -\ dL_+,
$$
and, when $\vert q_i\vert \ge 4$, the coefficients of $L_+$ 
(which have the special form of Proposition 4.1) may be found 
recursively by  performing finitely many integrations.
This algorithm is explained in \S 2 of \cite{Am-Gu}.

Let us now look at the two main families of examples in more detail (we
postpone comments on the case where $\deg q_i = 2$ to the end of
this section).

\no {\it 1. Full flag manifolds $G/B$}

For the algebra $\Cal A$ we take the deformation of the ordinary
cohomology algebra whose relations are the conserved quantities
of the open one-dimensional Toda lattice.  It may seem that
we are \ll starting with the answer\rrr, since this algebra has
already been identified with the quantum cohomology of $G/B$ in
\cite{Ki}, but our point of view here is that this algebra
exists naturally without reference to quantum cohomology.
We have $\vert q_i\vert=4$ for all $i$.

To construct the D-module $M^h$ we use the conserved quantities
of the  open one-dimensional {\it quantum} Toda lattice --- see
\cite{Ki} and \cite{Mr} for the precise definition.  These
are commuting differential operators which also have been studied
independently of quantum cohomology theory.  In particular,
it follows from \cite{Go-Wa} and the remarks at the beginning of this
section that $M^h$ is free over $K[h]$ with rank equal to $\dim H^\ast G/B$.  This is a quantization of $M^h$,
and so our method produces a \ll quantum cohomology algebra\rr and a
\ll quantum product operation\rrr.    Summarizing:

\proclaim{Theorem 4.3} The D-module $M^h$ associated to the
open one-dimensional
quantum Toda lattice is a quantization (in the sense of Definition 3.1)
of the algebra $\Cal A$ associated to the 
open one-dimensional Toda lattice. Hence  we obtain a
\ll quantum product\rr on $H^\ast G/B$ which may
be computed explicitly by the method explained above.
\qed
\endproclaim

Using the fact (\cite{Ki}) that $M^h$ is
known to be a quantization of the usual quantum cohomology algebra of $G/B$,
it can be shown (see \cite{Am-Gu}) that $L_+$ can be chosen to satisfy
$L_+\vert_{q=0}=I$, and furthermore that
our quantum product agrees with the
usual quantum product.
Computations for $G=\glnc$ $(n=2,3,4)$
are also given in \cite{Am-Gu}.

If the Schubert polynomial basis of $H^\ast G/B$ is used instead of the
monomial basis, then this procedure gives the so called quantum
Schubert polynomials.  Thus, $Q_0$ is essentially the \ll quantization map\rr
of \cite{FGP} and \cite{Ki-Ma} (for the case $G=\glnc$).
This theory has been well studied, but our approach 
makes clear why such a rich structure can be expected, and
in particular why the quantum products can be computed from surprisingly 
minimal assumptions about quantum cohomology. 

Finally, we should point out that the role of D-modules in
the approach of \cite{Gi-Ki},  \cite{Ki} to the computation of the
quantum cohomology algebra of $G/B$ (see also \cite{Mr}) is quite
different.
The main step there is to show that the conserved
quantities $D^h_i$ of the quantum Toda lattice imply relations 
$\lim_{h\to 0} S(D^h_i)$ of the quantum cohomology algebra
(in the notation of Definition 3.1). This uses the special fact that
the differential operators $D^h_i$ commute.

\no {\it 2. Fano toric manifolds with $\vert q_1\vert,\dots,\vert q_r\vert \ge 4$}

For the algebra $\Cal A$ we take the \ll provisional\rr quantum cohomology
algebra of Batyrev (\cite{Ba1}).  This exists for Fano and non-Fano toric
manifolds alike.  To construct a quantization we
shall use the theory of generalized hypergeometric partial differential equations of Gelfand,
Kapranov and Zelevinsky (\cite{GKZ},  \cite{HLY},  \cite{SST},  \cite{Co-Ka}). This
theory associates to a certain polytope a system of partial differential equations
or D-module, which we refer to as a GKZ D-module. Now, by a well known construction (see \cite{Od}),
such a polytope gives rise to a toric variety $M$ with a line bundle. 
We shall use this to prove:

\proclaim{Theorem 4.4} Let $M$ be a Fano toric manifold. Then there is a 
GKZ D-module which is
a quantization $M^h$ (in the sense of Definition 3.1) of Batyrev's algebra $\Cal A$.
Hence we obtain a \ll quantum product\rr on $H^\ast M$ which may be computed explicitly
by the method explained above.
\endproclaim

\demo{Proof} We need a GKZ D-module $M^{\sssize GKZ}$ whose rank is equal to the
dimension of the vector space $H^\ast M$.
The construction of suitable differential operators (defining  $M^h$) may then be
carried out exactly as in Section 5.5 of \cite{Co-Ka}, and it is easy to
see that these satisfy the conditions of Definition 3.1.

To obtain $M^{\sssize GKZ}$ we need a suitable polytope.  It is known
(see Lemma 2.20 of \cite{Od} and Section 2 of \cite{Ba2}) that, for
a Fano toric manifold, there exists a reflexive polytope
which gives rise to $M$ and has the following property: in the 
decomposition of the polytope given by taking the cones on the maximal faces
with common vertex at the origin, each such cone has unit volume. Therefore,
the volume of the polytope is the number of maximal faces, which (because
the polytope is reflexive) is equal to the number of maximal cones in a fan
defining the toric variety, and this in turn (by standard
theory of toric varieties) is equal to the number of fixed points of the
action of the torus on $M$.  This number is equal to the Euler characteristic of
$M$, and hence to $\dim H^\ast M$. On the other hand, it was
proved in \cite{GKZ},  \cite{SST} that the GKZ system in this situation
is free, with rank equal to the volume of the polytope.  We conclude that
the rank of $M^{\sssize GKZ}$ (and hence of $M^h$) is equal to $\dim H^\ast M$.
\qed\enddemo

To compute our quantum product explicitly, the
method of \cite{Am-Gu} may be used, exactly as in the case $M=G/B$.
To establish agreement with the
usual quantum product,
the method of \cite{Am-Gu} applies if one uses the fact that the that the GKZ D-module quantizes the
usual quantum cohomology algebra. 
This fact is known from very
general arguments (essentially, the mirror theorem
of Givental, as explained in Example 11.2.5.2 of \cite{Co-Ka}).
The simpler method used in \cite{Ki} in the case 
$M=G/B$ cannot be used in the Fano toric case, because the GKZ differential
operators do not in general commute.

We have assumed so far that $\vert q_i\vert\ge 4$ for all $i$.  If
some $\vert q_i\vert= 2$, the method of this section still applies, but
in Proposition 4.1 we have

\no(1) $Q_0=\exp X$ where $X_{\al,\be}=0$ for $\al\ge\be$,

\no(2) for $i\ge 1$, $(Q_i)_{\al,\be}=0$ for $\al\ge\be-i$.

\no In Corollary 4.2 we have 
$\hat c_0 = c_0 = 1$ and $c_i=b_i$ for $i=1,\dots,r$,
but $\hat c_i$ will in general be of the form $b_i+\sum a_j q_j$ 
(summing over $j$ such that $\vert q_j\vert= 2$).  
A similar phenomenon occurs for non-Fano manifolds, which are
the subject of the next section.

\head
\S 5 Beyond Fano manifolds
\endhead

Even for non-Fano manifolds, an algebra $\Cal A$ and a
quantization $M^h$ (with suitable coefficient algebra $K$)
lead to a gauge transformation $L_+=Q_0+O(h)$
and a connection $d+\hat\om$ with
$d\hat \om = \hat\om\wedge\hat\om = 0$.  
However, we do not necessarily have $\hat c_i=c_i=b_i$
for $i=1,\dots,r$, so we are not simply making a change of
basis in the algebra $\Cal A$.  We shall see that an important
new feature in the non-Fano case is the appearance of a coordinate
transformation (\ll mirror transformation\rrr).  

Referring to the proof of Proposition 3.2, let us define
$$
\align
\tilde L_-(q,h)&=Q_0(q)(I+h^{-1}A_1(q)+h^{-2}A_2(q)+\dots)\\
\tilde L_+(q,h)&= I + h Q_1(q)  + h^2 Q_2(q)  + \dots
\endalign
$$
i.e.\  we modify $L_-,L_+$ by moving the $Q_0$ factor from $L_+$ to $L_-$.  In
this case the proof shows that ${\tilde L_-}^{-1}d{\tilde L_-}$ is linear in $1/h$.  
Since the constant term of $\tilde L_+$ is the identity matrix, the
gauge transformation by ${\tilde L_+}$ simply changes the basis of $\Cal A$
as in \S 4.  In principle, therefore, it suffices to study the case
$$
\Om^h = \frac 1h\om + \th.
$$

Our first observation concerning this case is that 
the (usually complicated)  computation of $L_+$ becomes
easy.  Namely, we have $L_+=Q_0$ where $Q_0$ is a solution of
$Q_0^{-1}dQ_0 = \th$.
Then $\hat\Om^h=\frac1h\hat\om$ where 
$\hat\om = Q_0\om Q_0^{-1}$. The next step depends
on the form of $Q_0$. For example:

\proclaim{Proposition 5.1} Assume that 
$\vert q_1\vert,\dots,\vert q_r\vert \ge 0$ and that
$M^h=D^h/(D^h_1,\dots,D^h_u)$
is a quantization of $\Cal A$, where each $D^h_i$ is 
homogeneous in the sense of Proposition 3.3.
Assume further that $\Om^h = \frac 1h\om + \th$ and that
the zero-th order term of any second order element of a Gr\"obner basis of
$(D^h_1,\dots,D^h_u)$ is independent of $h$.
Then the block structure of $Q_0$ has the form
$$
Q_0=
\pmatrix
1 & 0 & [\ast] \\
0 & T & [\ast] \\
[0] & [0] & [\ast]
\endpmatrix
$$
where  $[\ast]$ denotes a submatrix and
$[0]$ denotes a zero submatrix (where a submatrix may
consist of several blocks).  That is,
$(Q_0)_{\al,\be}=0$ if $\al=0$ or $1$, unless $(\al,\be)=(0,0)$ or $(\al,\be)=(1,1)$.
\endproclaim

\demo{Proof} Since $Q_0^{-1}dQ_0=\th$, it suffices to prove that $\th$ has the
block form
$$
\pmatrix
0 & 0 & [\ast] \\
0 & \ast & [\ast] \\
[0] & [0] & [\ast]
\endpmatrix,
$$
since the matrices of this type form a Lie algebra.   From the definition
of $\Om^h$, the nonzero entries of $\th$ arise from expressions of
the form $h\b_i P_j$ which contain terms with \ll excess $h$\rrr, i.e.\  terms
which still contain $h$ after replacing $h\b_1,\dots,h\b_r$ by $b_1,\dots,b_r$.
Since $P_0=1$ we have $h\b_i P_0 =h\b_i$, and by the definition of
quantization (cf.\  the proof of Corollary 4.2) the
reduction of $h\b_i$ modulo these generators is $h\b_i$ itself.  
There are no excess $h$ here, so the first
column of $\th$ is zero.  Regarding the second column of $\th$, the
third sub-matrix is zero by (2a) of Proposition 3.3:  each block is homogeneous
of negative degree and well defined at $q=0$, hence polynomial in
$q_1,\dots,q_r$; but this contradicts the assumption
$\vert q_1\vert,\dots,\vert q_r\vert \ge 0$, unless that block is
zero.  By assumption,  there are no
excess $h$ in the zero-th order term of $h^2\b_i \b_j$, so the
first block is also zero, as required.
\qed\enddemo

Our second observation is that, while
$T$ is not necessarily the identity matrix, it does have a special form:

\proclaim{Proposition 5.2} Assume that $M^h=D^h/(D^h_1,\dots,D^h_u)$
is as in Proposition 5.1.  Then the matrix $T$ is a Jacobian matrix,
i.e.\  there exist new local coordinates $\hat t_1,\dots,\hat t_r$ on
the vector space $\C^r$ such that
$$
T=
\pmatrix
\b_1\hat t_1& \cdots & \b_r\hat t_1 \\
\vdots & & \vdots \\
\b_1\hat t_r& \cdots & \b_r\hat t_r
\endpmatrix.
$$ 
\endproclaim

\demo{Proof} 
By definition, $\b_i P_j = \sum_{k=0}^s (\Om^h_i)_{kj} P_k$
mod $(D^h_1,\dots,D^h_u)$.
We have $h\b_i P_j = h^2\b_i\b_j$ for $1\le i,j\le r$.  Since
$\b_i\b_j=\b_j\b_i$, it follows that $(\Om^h_i)_{kj} = (\Om^h_j)_{ki}$
for $1\le i,j\le r$. In particular this symmetry is valid for $\th$,
and for the $(1,1)$ block of $\th=Q_0^{-1}dQ_0$, namely for $T^{-1}dT$.  It is easy to
verify that this implies that the operators 
$\hat\b_i = \sum_{j=1}^r  (T^{-1})_{ji}\b_j$
($i=1,\dots,r$) commute, hence define new local coordinates.
\qed\enddemo

Under the coordinate transformation $t\mapsto \hat t$, the differential
operators $D^h_i$ transform to differential operators $\hat D^h_i$.
Let $\hat D^h$ be the ring of differential operators analogous to $D^h$,
using $\hat\b_i= \b/\b \hat t_i = \hat q_i \b/\b \hat q_i$ instead of $\b_i$.
Then
we obtain a new D-module $\hat M^h= \hat D^h/(\hat D_1^h,\dots,\hat D_u^h)$
and a de-quantized commutative algebra $\hat \Cal A$.  With respect
to the basis of standard monomials in 
$h\hat\b_1,\dots,h\hat\b_r$ we obtain a
connection $\hat d+\tilde\Om^h$, and by construction $\tilde\Om^h$
has the property $\tilde T = I$.
We can now apply the procedure of \S 3 to $\hat M^h$.
A gauge transformation  produces a connection 
$\hat d+\hat\Om^h$  with 
$\hat\Om^h=\frac 1h \hat\om$, so we can define a
\ll quantum product operation\rr on $\Cal A_0$.

This is our general procedure for reconstructing quantum cohomology:
first we make a change of variable to obtain a connection of the kind
discussed in \S 3, then we make a gauge transformation to
obtain a connection with the properties of the Dubrovin connection.
The first operation is
natural from the point of view of the D-module $M^h$, but it does not in
general preserve the isomorphism type of the associated algebra $\Cal A$.
The second one preserves this isomorphism type, and just introduces the
additional information needed to define a quantum product.  We summarize
this in the following theorem.

\proclaim{Theorem 5.3}  Assume that 
$M^h= D^h/(D_1^h,\dots,D_u^h)$ is a quantization of $\Cal A$, 
that the conditions of Proposition 5.1 hold, and that $\Om^h = \frac 1h\om + \th$.
Then by a change of variable and a gauge transformation we obtain a connection form
$\hat\Om^h=\frac 1h \hat\om$ satisfying 
$\hat d \,\hat\om = \hat\om \wedge \hat\om = 0$,
and a \ll quantum product operation\rr on $\Cal A_0 = H^\ast M$.
\qed
\endproclaim

\no{\it Example 5.4:} The Hirzebruch surfaces 
$\Si_k = \P(\Cal O(0)\oplus \Cal O(-k) )$, where $\Cal O(i)$
denotes the holomorphic line bundle on $\C P^1$ 
with first Chern class $i$, are Fano when $k=0,1$. 
We shall consider the first non-Fano
case, $\Si_2$.  The ordinary cohomology algebra is
$$
\Cal A_0 = H^\ast \Si_2 = \C[b_1,b_2]/(b_1^2,b_2(b_2-2b_1))
$$
(in the notation of \cite{Gu2}, $b_1=x_1$ and $b_2=x_4$.)  This is
a complex vector space of dimension $4$.
Batyrev's algebra (\cite{Ba1}), obtained by consideration
of rational curves in $\Si_2$,  is in this case
$$
\Cal A=\C[b_1,b_2,q_1,q_2]/(b_1^2-q_1(b_2-2b_1)^2,b_2(b_2-2b_1)-q_2).
$$
It is a $\C[q_1,q_2]$-module of rank $4$.  We have $\vert q_1\vert = 0$ 
and $\vert q_2\vert = 4$ here.

Consider the D-module $M^h=D^h/(D^h_1,D^h_2)$ where
$$
D^h_1 = h^2\b_1^2 - q_1h^2(\b_2-2\b_1)(\b_2-2\b_1-1),
\quad
D^h_2 = h^2\b_2(\b_2-2\b_1) - q_2.
$$
This can be derived from a GKZ D-module, as in the Fano case
(see \cite{Co-Ka}, Section 5.5).
It is a  D-module which is free over $K[h]$ of rank $4$,
where $K$ is the field of rational functions, and therefore a quantization of $\Cal A$.
(It is interesting to note that the \ll naive quantization\rrr, obtained by
using $D^h_1 = h^2\b_1^2 - q_1h^2(\b_2-2\b_1)^2$
and $D^h_2 = h^2\b_2(\b_2-2\b_1) - q_2$, has rank $0$, and
is therefore {\it not} a valid quantization of $\Cal A$.)  The 
Gr\"obner basis for the ideal $(D^h_1,D^h_2)$,  turns out to be
{
$$
\gather
\underline{2 h^2 \b_1\b_2} -  h^2 \b_2^2  + q_2,\\
\underline{(4q_1-1) h^2\b_1^2} - q_1h^2 \b_2^2 + 2q_1 h^2 \b_1 - q_1 h^2 \b_2 +2q_1q_2,\\
\underline{h^3\b_2^3} + 2q_2(4q_1-1) h\b_1 - q_2(4q_1+1) h\b_2 - hq_2
\endgather
$$
}
The equivalence classes of the standard monomials
$
1,h\b_2,h\b_1,h^2\b_2^2
$
(i.e.\  the monomials  $(h\b_1)^i(h\b_2)^j$ not \ll divisible\rr by any of the
leading terms, which are underlined)
form a basis of $M^h$. With respect to this basis, the matrices
$\Om^h_i$ (of the action of $\b_i$) are:
$$
{
\Om^h_1 =
\frac1h
\pmatrix
0 & -\frac{q_2}2 & \frac{-2q_1q_2}{4q_1-1}
& 0 \\
0 & 0 & h \frac{q_1}{4q_1-1} & 2q_1q_2 \\
1 & 0 & h \frac{-2q_1}{4q_1-1} & -q_2(4q_1-1) \\
0 & \frac12 & \frac{q_1}{4q_1-1} & 0
\endpmatrix,
\ \ 
\Om^h_2 =
\frac1h
\pmatrix
0 & 0 & \frac{-q_2}{2} & h q_2 \\
1 & 0 & 0 & q_2(4q_1+1) \\
0 & 0 & 0 & -2q_2(4q_1-1) \\
0 & 1 & \frac12 & 0
\endpmatrix}
$$
In particular we see that $\Om^h$ is of the form $\frac1h\om+\th$ here. 

The gauge transformation $L_+=Q_0$ such that $Q_0^{-1}dQ_0=\th$ and 
$Q_0\vert_{q=0}=I$ is easily found.  Its inverse is
$$
{
Q_0^{-1}=
\pmatrix
1 & 0 & 0 & -q_2 \\
0 & 1 & \frac12(1-\sqrt{1-4q_1}) & 0\\
0 & 0 & \sqrt{1-4q_1} & 0 \\
0 & 0 & 0 & 1
\endpmatrix}
$$
The coordinate transformation is determined by the 
central $2\times 2$ block of $Q_0^{-1}$, i.e.\ 
$$
\hat\b_2 = \b_2,
\quad
\hat \b_1 = \frac12(1-\sqrt{1-4q_1})\b_2 +  \sqrt{1-4q_1} \b_1.
$$
Writing $\hat q_i = e^{\hat t_i}$, it is easy to deduce that
$q_1={\hat q_1}/{(1+\hat q_1)^2}$, $q_2 = \hat q_2(1+\hat q_1)$,
if we impose the condition that the origin maps to the origin.
Let us see what effect this transformation has on the D-module $M^h$. 
From
$$
\b_2 =\hat\b_2,
\quad
\b_1 =  - \frac{\hat q_1}{1-\hat q_1}\hat\b_2 +
\frac{1+\hat q_1}{1-\hat q_1}\hat\b_1
$$
we obtain
$$
\align
D^h_1 &= h^2\hat\b_1^2 - \hat q_1\hat q_2
-
\frac{\hat q_1}{1-\hat q_1}
\big(
h^2\hat\b_2 (\hat \b_2 - 2\hat\b_1)  -  \hat q_2(1-\hat q_1)
\big)
\ (=\hat D^h_1,\ \text{by definition})
\\
D^h_2 &=
\frac{1+\hat q_1}{1-\hat q_1}
\big(
h^2\hat\b_2 (\hat \b_2 - 2\hat\b_1)  -  \hat q_2(1-\hat q_1)
\big) \ (=\hat D^h_2,\ \text{by definition}).
\endalign
$$
These operators define an equivalent D-module $\hat M^h$, but the
de-quantized algebra $\hat\Cal A$ is quite different from $\Cal A$:
$$
\hat\Cal A=
\C[\hat b_1,\hat b_2,\hat q_1,\hat q_2]/
(\hat b_1^2-\hat q_1\hat q_2,
\hat b_2 (\hat b_2-2\hat b_1)-\hat q_2(1-\hat q_1)).
$$
To obtain a \ll quantum product operation\rr on $H^\ast \Si_2$
we carry out the procedure of \S 3, 
but this time starting
from $\hat\Cal A$.  With respect to the standard monomials 
$1,h\hat\b_2,h\hat\b_1,h^2\hat\b_2^2$, the matrices $\tilde\Om^h_i$
(of the action of $\hat\b_i$) can be computed as
$$
{
\tilde\Om^h_1 =
\frac1h
\pmatrix
0 & -\frac{\hat q_2(1-\hat q_1)}2 & \hat q_1\hat q_2
& h \hat q_1\hat q_2\\
0 & 0 & 0 & 2\hat q_1\hat q_2 \\
1 & 0 & 0 & \hat q_2(1-\hat q_1) \\
0 & \frac12 & 0 & 0
\endpmatrix,
\quad
\tilde\Om^h_2 =
\frac1h
\pmatrix
0 & 0 & -\frac{\hat q_2(1-\hat q_1)}2 & h \hat q_2(1+\hat q_1) \\
1 & 0 & 0 & \hat q_2(1+3\hat q_1) \\
0 & 0 & 0 & 2\hat q_2(1-\hat q_1) \\
0 & 1 & \frac12 & 0
\endpmatrix}
$$
We have $\tilde\Om^h=\frac1h\tilde\om+\tilde\th$.
The inverse of the matrix $\tilde Q_0$ such that $\tilde Q_0^{-1}d\tilde Q_0=\tilde\th$ and 
$\tilde Q_0\vert_{\hat q=0}=I$ is
$$
{
\tilde Q_0^{-1}=
\pmatrix
1 & 0 & 0 & -\hat q_2(1+\hat q_1) \\
0 & 1 & 0 & 0\\
0 & 0 & 1 & 0 \\
0 & 0 & 0 & 1
\endpmatrix}
$$
This converts $\tilde \Om^h$ to $\hat \Om^h = \frac 1h \hat\om$, where
$\hat\om=\tilde Q_0\tilde\om \tilde Q_0^{-1}$, and we have:
$$
{
\hat\om_1=
\pmatrix
0 & \hat q_1\hat q_2 & \hat q_1\hat q_2&0\\
0 &0 &0 &2\hat q_1\hat q_2\\
1 &0 & 0& -2\hat q_1\hat q_2\\
0 &\frac12 &0 &0
\endpmatrix,
\quad
\hat\om_2=
\pmatrix
0 & (1+\hat q_1)\hat q_2 & \hat q_1\hat q_2&0\\
1 &0 &0 &2\hat q_1\hat q_2\\
0 &0 & 0& -2\hat q_2(-1+\hat q_1)\\
0 &1 &\frac12 &0
\endpmatrix.}
$$
We obtain the following basic products:   
$ \hat b_1 \circ_{\hat t}  \hat b_1  = \hat q_1\hat q_2$,  
$\hat b_1 \circ_{\hat t}  \hat b_2  =  \hat b_1\hat b_2  +\hat q_1\hat q_2 $, 
$\hat b_2 \circ_{\hat t}  \hat b_2  =  \hat b_2^2 +  \hat q_2(1+\hat q_1) $.
These are in agreement with the observation made at the end of Chapter 11
of \cite{Co-Ka} that the quantum products of $\Si_2$ can be deduced
from those of $\Si_0=\C P^1\times\C P^1$, if one uses the symplectic
invariance of Gromov-Witten invariants.  Thus our product is indeed the
usual quantum product.
\qed

The coordinate transformation (\ll mirror transformation\rrr) in this example was
obtained in Example 11.2.5.2 of \cite{Co-Ka}, as a consequence of
Givental's \ll Toric Mirror Theorem\rrr.
It appeared originally, in a similar situation, in the Introduction to \cite{Gi4}.
In fact, as we shall discuss elsewhere,
this example is typical of the case of a semi-positive toric manifold
$M$, i.e.\  a toric manifold such that the evaluation of
$c_1M$ on any homology class represented by a rational curve is non-negative.

There are two ways to apply our theory in this situation. 
The first point of view is to assume that the quantum cohomology of $M$
is given by Givental's mirror theorem, then use this to deduce that our quantum
cohomology agrees with the usual quantum cohomology.  Alternatively, our
construction of quantum cohomology can be used to prove a version of
the mirror theorem.

To explain the latter, we need the
explicit solution $J_{\sssize GKZ}$ of the GKZ system constructed 
in \cite{GKZ}, \cite{St},  \cite{SST} (the function
$I_{\nu}$ in (11.73) of \cite{Co-Ka} with $\nu=0$ and $t_0=0$).
Let 
$$
L_{\sssize GKZ}=
\pmatrix
\vert &  & \vert \\
P_0 J_{\sssize GKZ} & \cdots &  P_s J_{\sssize GKZ} \\
\vert &  & \vert
\endpmatrix.
$$
This is the (transpose of the) fundamental solution of the first order system
$(d-(\Om^h)^t)L^t=0$. (Since $P_0=1$
we have $P_0 J_{\sssize GKZ} =  J_{\sssize GKZ}$, of course.)  For a
semi-positive toric manifold we may apply the method of this section to
$L=L_{\sssize GKZ}$.  We obtain a new first order system
$(\hat d-(\hat\Om^h)^t)\hat L^t=0$, with fundamental solution of the form
$$
\hat L=
\pmatrix
\vert &  &  \\
\hat J & \cdots &   \\
\vert &  & 
\endpmatrix.
$$
The relation between $\hat\Om^h$ and $\Om^h$ is
$\hat\Om^h=G^{-1}(X\Om^h)G - dG G^{-1}$, where $G$ is a gauge transformation
and $X$ is the matrix function expressing the relation between the standard
monomial bases of $\hat M^h$ and $M^h$.  From our earlier descriptions of $G$ and $X$,
it is obvious that the first rows of $\hat\Om^h, \Om^h$ (i.e.\  the first columns of 
$(\hat\Om^h)^t, (\Om^h)^t$) are unaffected by $G$ or $X$. Hence 
$$
\hat J(\hat t\,) =  J_{\sssize GKZ}(t).
$$
This is (our version of) Givental's toric mirror theorem.  It expresses a relation between
the structure constants (\ll Gromov-Witten invariants\rrr) for our quantum product
operation and the coefficients of the generalized hypergeometric series 
$ J_{\sssize GKZ}$.  Of course this is merely a reflection of the more fundamental underlying
relation between $\hat M^h$ and $M^h$.

\bigskip
\it
\no  
Department of Mathematics,
Graduate School of Science,
Tokyo Metropolitan University,
Minami-Ohsawa 1-1, Hachioji-shi,
Tokyo 192-0397, Japan

\medskip
\no 
martin\@comp.metro-u.ac.jp

\enddocument